\documentstyle[11pt,amsfonts]{article}

\begin{document}

\bibliographystyle{plain}

\newcommand{\comment}[1]{}

\def\trdeg{\mbox{\normalshape tr.deg}}
\def\c{{\Bbb C}}
\def\R{{\Bbb R}}
\def\Rexp{{\Bbb R}_{\exp}}
\def\q{{\Bbb Q}}
\def\n{{\Bbb N}}
\def\z{{\Bbb Z}}

\newtheorem{conjecture}{Conjecture}
\newtheorem{theorem}{Theorem}
\newtheorem{lemma}{Lemma}
\newtheorem{proposition}{Proposition}
\newtheorem{corollary}{Corollary}
\newtheorem{definition}{Definition}
\newtheorem{problem}{Problem}
\newtheorem{remark}{Remark}
\newtheorem{example}{Example}
\newtheorem{hypothesis}{Hypothesis}

\makeatletter
\def\@yproof[#1]{\@proof{ #1}}
\def\@proof#1{\begin{trivlist}\item[]{\em Proof#1.}}
\newenvironment{proof}{\@ifnextchar[{\@yproof}{\@proof{} 
}}{~$\Box$\end{trivlist}}
\makeatother

\def\k{{\Bbb K}}


\title{On Defining Irreducibility}

\author{Pascal Koiran\\
Laboratoire de l'Informatique du Parall\'elisme\\
Ecole Normale Sup\'erieure de Lyon -- CNRS\\
46, all\'ee d'Italie\\
69364 Lyon Cedex 07, France\\
\tt Pascal.Koiran@ens-lyon.fr}

\maketitle

In this note we prove the following result.
\begin{theorem} \label{main}
For any $n \geq 2$, irreducibility is not a definable property of 
real algebraic varieties of $\R^n$.
\end{theorem}
This means that for all $n \geq 2$, there exists no
first-order sentence 
$F_n$ in the language ${\cal L}_n = (+,-,\times,=,0,1,I_n)$
satisfying the following property:
for any real algebraic variety $V$ of $\R^n$, if we interpret the
$n$-ary predicate $I_n$ by membership to $V$ then $\R \models F_n$
if and only if $V$ is irreducible. Our proof will imply that
irreducibility remains undefinable even if we add the exponential
function to ${\cal L}_n$. It was inspired by Wilkie's proof that $C^{\infty}$
functions 
are not definable~\cite{Wilkie94}.
\begin{conjecture} 
For any $n \geq 2$, irreducibility is not a definable property of 
complex algebraic varieties of $\c^n$.
\end{conjecture}
It is not clear how this conjecture should be attacked. 
In particular, if one wishes to follow the same proof strategy as in
this note, it is not clear by what structure the real exponential
field should be replaced.

A study of definability in algebraically closed fields was initiated 
in~\cite{ChaKoi98}: in that paper we gave some examples of definable properties
(for instance, dimension is definable), showed that connectivity of 
algebraic varieties is not definable, and conjectured that algebraic
varieties are not definable among constructible sets.
There is already a fairly large body of work on definability 
in o-minimal and stable structures, 
see for instance~\cite{BaldBe98,BST99,BDLW98,BeLib96,BeLib97,GS97}.

Let $V_n$ be the real variety $\{(x,y) \in \R^2;\ P_n(x,y)=0\}$
where $$P_n(x,y) = x^{2n} - (1+y^2)^2.$$ Note that 
$P_n =Q_n R_n$ with $Q_n(x,y)=x^n-1-y^2$ and $R_n(x,y)=x^n+1+y^2$.
Note also that there is a formula $F(n,x,y)$ of  
the language ${\cal L}_{\exp}=\{+,-,\times,=,\exp,0,1\}$ such that
whenever $n$ is a positive integer, $F(n,.,.)$ defines $V_n$:
take for $F$ the formula
$$\exists z [e^z=x^2 \wedge e^{nz} = (1+y^2)^2].$$
Here is another useful remark.
\begin{lemma} \label{Qn}
For any integer $n \geq 0$, the real variety $V_{\R}(Q_n)$ 
defined by $Q_n$ is irreducible.
\end{lemma}
\begin{proof}
For $n=0$, $V_{\R}(Q_n)$ is the line $\{y=0\}$ and is therefore irreducible.
For $n \geq 1$, $V_{\R}(Q_n)$ also has dimension 1 and it is not hard
to check that $Q_n$ is irreducible over $\c$. It therefore follows from 
Lemma~\ref{irred} below that $V_{\R}(Q_n)$ is irreducible.
\end{proof}
\begin{lemma} \label{irred}
Let $Q \in \R[x,y]$ be a polynomial defining a real variety
$V_{\R}(Q)$ of dimension 1. 
If $Q$ is irreducible over $\c$ then $V_{\R}(Q)$ is irreducible.
\end{lemma}
\begin{proof}
Let $V_1$ and $V_2$ be two real varieties such that 
$V_{\R}(Q)=V_1 \cup V_2$. We shall see that 
$V_1 =V_{\R}(Q)$ or $V_2 =V_{\R}(Q)$. 
At least one of the $V_i$ must be of
dimension 1. Assume for instance that this is true of $V_1$,
and let $P_1 \in \R[x,y]$ be a polynomial such that $V_1=V_{\R}(P_1)$.
Since $V_{\R} (P_1) \cap V_{\R}(Q)$ has dimension 1, 
the complex variety $V_{\c} (P_1) \cap V_{\c}(Q)$ also has dimension 1. 
By irreducibility of $Q$, $V_{\c}(Q)$ is irreducible. 
It follows that  $V_{\c}(Q) \subseteq V_{\c} (P_1)$ 
and therefore $V_{\R}(Q) \subseteq V_{\R}(P_1)$, 
i.e., $V_{\R}(P_1)=V_{\R}(Q)$.
\end{proof}
\begin{proposition} \label{parity}
$V_n$ is irreducible if and only if $n$ is even.
\end{proposition}
\begin{proof}
If $n$ is an odd number, the real varieties defined by $Q_n$ and $R_n$
 are distinct (they are in fact disjoint).
Hence $V_n$ is not irreducible in this case.

If $n$ is even, $R_n$ has no real zeros. $V_n$ is therefore
irreducible by Lemma~\ref{Qn}.
\end{proof}

\begin{proof}[of Theorem~\ref{main}]
We shall prove that irreducibility is not definable in dimension 2.
The result for higher dimensions follows immediately from this special
case.

Assume by contradiction that there exists a formula $F_2$ 
of ${\cal L}_2$ such that for any real variety $V$ of $\R^2$, 
if we interpret $I_2$ by membership to $V$ then $\R \models F_2$
if and only if $V$ is irreducible.
Let $G(n)$ be the formula obtained from $F_2$ by replacing each
instance $I_2(x,y)$ of $I_2$ in $F_2$ by $F(n,x,y)$.
By Proposition~\ref{parity}, $G(n)$ is true whenever $n$ is an even
integer, and is false whenever $G(n)$ is an odd integer.
This is in contradiction with the 
o-minimality of the real exponential field (see~\cite{Wilkie98} 
and the references there).
\end{proof}

\end{document}